\documentstyle[amsfonts,twoside,11pt]{article}

\newtheorem {theo} {\bf Theorem} [section]
\newtheorem {prop} [theo] {\bf Proposition}
\newtheorem {cory} [theo] {\bf Corollary}
\newtheorem {lem} [theo] {\bf Lemma}
\newtheorem {defn} [theo] {\bf Definition}

\newcommand{\QED}{\hfill \CaixaPreta \vspace{6mm}}
\def\CaixaPreta{\vrule Depth0pt height6pt width6pt}

\newcommand{\dpy}{\displaystyle}
\newcommand{\be}{\begin{eqnarray}}
\newcommand{\ee}{\end{eqnarray}}
\newcommand{\benn}{\begin{eqnarray*}}
\newcommand{\eenn}{\end{eqnarray*}}
\newcommand{\bse}{\begin{equation}}
\newcommand{\ese}{\end{equation}}
\newcommand{\bsenn}{\begin{displaymath}}
\newcommand{\esenn}{\end{displaymath}}

\newcommand{\logand}{\;\;{\rm and }\;\;}
\newcommand{\Hd}{\;{\rm Hd }\,}
\newcommand{\Hdim}{\;{\rm Hdim }\,}

\newcommand{\mod}{\bmod}
\newcommand{\such}{\;\;{\rm such\; that }\;\;}
\newcommand{\where}{\;\;{\rm where }\;\;}
\newcommand{\Id}{\,{\rm Id}\,}

\renewcommand{\O}{\mbox{$\Bbb O$}} 	
\newcommand{\G}{\mbox{$\Bbb G$}} 	
\newcommand{\R}{\mbox{$\Bbb R$}} 	
\renewcommand{\S}{\mbox{$\Bbb S$}}  
\def\SBIMSMark#1#2#3{
 \font\SBF=cmss10 at 10 true pt
 \font\SBI=cmssi10 at 10 true pt
 \setbox0=\hbox{\SBF Stony Brook IMS Preprint \##1}
 \setbox2=\hbox to \wd0{\hfil \SBI #2}
 \setbox4=\hbox to \wd0{\hfil \SBI #3}
 \setbox6=\hbox to \wd0{\hss
             \vbox{\hsize=\wd0 \parskip=0pt \baselineskip=10 true pt
                   \copy0 \break%
                   \copy2 \break%
                   \copy4 \break}}
 \dimen0=\ht6   \advance\dimen0 by .75\vsize \advance\dimen0 by 8 true pt
 \dimen2=\hsize \advance\dimen2 by .25 true in
 \ht6=0pt \dp6=0pt
 \setbox8=\vbox to \dimen0{\vfill \hbox to \dimen2{\hss \copy6}}
 \ht8=0pt \dp8=0pt \wd8=0pt
 \copy8
}

\begin{document}
\title{RIGIDITY PROPERTIES OF LOCALLY SCALING FRACTALS}
\author {J. J. P. Veerman\thanks{Rockefeller University, New York. e-mail:
veerman@golem.rockefeller.edu} \ and \ L. B. Jonker\thanks{Queen's University, Kingston. e-mail: leo@mast.queensu.ca}} 
\date{\relax}
\maketitle
\SBIMSMark{1997/2}{January 1997}{}

\vskip .2in

\begin{abstract}
A set has local scaling  if in a neighborhood of a point the structure
of the set can be mapped onto a finer scale structure of the set.
These scaling transformations are compact sets of locally affine 
 contractions (that is, contractions with uniformly $\alpha$-H\"older continuous derivatives). In this setting, without the open set condition or any other assumption on the spacing of these contractions, we
show that the measure of the set is an upper semi-continuous function
of the scaling transformations in the $C^0$-topology. With a restriction on 
the 'non-conformality' of the scaling transformations  we show that the  Hausdorff dimension is a lower semi-continuous function in the $C^{1}$-topology.
We include some examples to show that continuity does not follow in either case.
\end{abstract} 

\vskip .3in

\section{Introduction}
\label{introduction}
\setcounter{figure}{0}
\setcounter{equation}{0}

\vskip.2in

We consider a certain class of contractions on the 
space $H(\bar{I})$ of compact subsets of the closed unit ball $\bar{I}$ in $\R^n$. 
For the purpose of this introduction, it is convenient
to think of such a contraction as a finite collection of contracting 
diffeomorphisms as in the examples at the end of this section. More general and precise  
definitions will be given in section \ref{definitions}.  For the moment the following definition of the contraction 
on the space of compact sets will suffice:

\begin{defn}
\label{def=ifs}
Let $\{f_j\}_{j\in J}$ be a finite collection of contracting 
homeomorphisms.  The associated  iterated function 
system is the map ${\cal F} : H(\bar{I}) \rightarrow H(\bar{I})$, determined by 
\begin{itemize}
\item ${\cal F}(x) = \cup_{j \in J} f_j(x)$ .
\item $A \in H(\bar{I})$ then ${\cal F}(A) = \cup_{x \in A} {\cal F}(x)$ .
\end{itemize}
\end{defn}

It will be convenient, and it will not lead to confusion, to use the symbol 
\,${\cal F}$\, for the collection of functions as well as for the iterated function system to which it gives rise.  We will do so frequently throughout the paper.

The fact that such iterated function systems are contractions was observed by \cite{Hu} and 
has a particularly important consequence:
Let $A \in H({\bar I})$, then $\{ {\cal F}^n (A)\}_n$ converges 
uniformly to a uniquely defined  compact set $\Lambda ({\cal F})$.

The object of our study is the  fixed point  $\Lambda ({\cal F})$ of these maps. By definition, 
then, this fixed point is a compact set for which we have local
scaling as expressed by the fixed point equation:
\bsenn
{\cal F}(\Lambda) = \Lambda \quad .
\esenn
Our interest centers on the question of how the measure and the Hausdorff dimension 
of this set vary as we vary the contraction ${\cal F}$.

To this end, we define  spaces $\G^0$ and $\G^1$ of such contractions and put topologies on them.   The space $\G^0$ is intended for a study of the Lebesgue measure of $\Lambda$, and is given a topology induced by the $C^0$-topology on the space of homeomorphisms.  $\G^1$ is the subspace consisting of iterated function systems generated by contracting diffeomorphisms with uniform H\"{o}lder continuous derivatives, and is intended  for  considerations
related to the Hausdorff dimension of $\Lambda$. We define two functions:

\begin{defn}
\label{def=measure}
Let $\mu : \G^0 \rightarrow \R^+$ be the function that assigns the Lebesgue measure 
of $\Lambda({\cal F})$ to ${\cal F}$. 
\end{defn}

\begin{defn}
\label{def=dimension-function}
Let $\Hdim : \G^1 \rightarrow \R^+$ be the function that assigns the Hausdorff dimension 
of $\Lambda({\cal F})$ to ${\cal F}$.
\end{defn}

In section \ref{definitions}, we define our notions (in particular the spaces
of contractions on $H(\bar{I})$) carefully. The space we
define to deal with the measure theoretic question is in fact much more 
general than here indicated. In this setting we then prove in section \ref{measure}
the following result.

\begin{theo}
The function $\mu$ is upper semi-continuous.
\end{theo}

The derivatives of the functions $\{f_j\}_{j\in J}$ that constitute the map 
${\cal F}$ are of vital importance for dimension 
calculations. So we must assume that the functions generating $\G^1$ have H\"{o}lder continuous derivatives. 
This space  is also defined
in section \ref{definitions}. 

As it happens, the theory involving derivatives is much subtler, and in addition requires  the assumption that the system is semi-conformal.

\begin{defn}
\label{def=semi-conformal} 
The subset $\S \subset \G^1$ is the collection of semi-conformal systems ${\cal F}$.  
A system ${\cal F} \in \G^1$ is called semi-conformal if 
\bsenn
\lim_{n\rightarrow \infty} \max_{f_i \in {\cal F}, \, x_0 \in {\bar I}}\dpy\frac{1}{n} \ln \parallel (Df_n \cdots f_1|_{x_0})\parallel \cdot \parallel (Df_n \cdots f_1|_{x_0})^{-1}\parallel  = 0 \quad . 
\esenn
\end{defn}

Note that if the derivatives of the functions are constant and equal,
then the assumption of semi-conformality is equivalent to the assumption that the eigenvalues are equal in modulus.  The assumption of semi-conformality is closely related to the main asumption in the definition of an asymptotic Moran symbolic construction in \cite{PW}.

The proof of the main theorem passes through estimates 
of ratios of derivatives of long compositions taken at different points.
These distortion calculations have been done before only in the case 
of one-dimensional systems. In section \ref{distortion}, we do these 
calculations for arbitrary dimension.

\begin{defn}
\label{def=open-set}
The subset $\O \subset \G^1$ is the collection of systems ${\cal F}$ satisfying the open 
set condition. That is, ${\cal F}$ consists of finitely many functions $f_i$
and there is an open set $V$ containing 
$\Lambda ({\cal F})$ such that $f_i(V) \cap f_j(V)$ 
is empty whenever $i\ne j$.
\end{defn}

In section \ref{dimension}, we prove the following:

\begin{theo} 
\label{thm=main2}
The function $\Hdim$ is continuous on $\S \cap \O \cap \G^1$ and lower semi-continuous on
$\S\cap \G^1$. 
\end{theo}

\begin{cory}
\label{thm=cap}
On $\S\cap\G^1$ the Hausdorff dimension and the limit capacity are equal.
\end{cory}

The definition of limit capacity is given at the end of this section.  It is an important concept because it lends itself to numerical calculations and because it is the concept used in embedding theorems (see \cite{SYC}).

The principal result here is the semi-continuity. The examples below will 
show that continuity does not hold. 
We note that in previous works (\cite{MM}, \cite{PV}, and \cite{Ta}), 
continuity of the dimension has been proved
if the measure function $\mu$ is restricted to conformal systems and if ${\cal F}$  
satisfies a certain condition similar to the open set condition. In these cases
smoothness of the dimension can also be proved (see for example \cite{Ru1}). 
Other workers have related the Hausdorff dimension to other quantities,
such as the Lyapunov exponents \cite{Yo}. In the absence of the open set 
condition, it is not clear how to define Lyapunov exponents. 
 
Here are some examples to illustrate the subtlety of the
problem we are studying. Note that they establish that the functions mentioned
 are not continuous.

For $t \in [0,1/2]$, let ${\cal F}_t$ be given by $\{f_i\}_{i=0}^2$ where
\bsenn
\left\{ \begin{array}{lll}
f_0(x) &=& \dpy\frac{x}{3} \\[.4cm]
f_1(x) &=& \dpy\frac{x+t}{3} \\[.4cm]
f_2(x) &=& \dpy\frac{x+1}{3}
\end{array} \right.
\esenn
Note that each function maps the unit interval into itself, thus the system
is a contraction of $H([0,1])$ into itself.
 
The dimension and measure of the invariant now  depend only on the parameter $t$.
 
\begin{theo}
\label{thm=3branch}
Let ${\cal F}_t$  be the system just described. Then \\
i) If $t=p/q$ is rational and $pq \mod 3 = 2$ then $\mu(t) = 1/q$.\\
ii) If $t=p/q$ is rational and $pq \mod 3 \neq 2$ then $\Hdim(t) < 1$.\\
iii) For all irrational $t$, $\mu(t) = 0$.\\
iv) For almost all $t$, $\Hdim(t) = 1$.
\end{theo}
 
We remark that in the rational case ii), dimension calculations, although
apparently feasible, are unknown to us.
 
The first result is a special case of a result proved in \cite{GH}, with a more
geometrical proof in \cite{Ve1}. The second statement is implied by Theorem 4.1
of \cite{LW} and Theorem 2.3 of \cite{Ve2}.
The third part of the theorem is due to Odlyzko \cite{Od} for
almost all $t$. In \cite{LW2} this was generalized to include all
irrational $t$. The last part is a special case of Theorem 9.12 of \cite{Fa}.

Another example is given by the following family of systems ${\cal G}_\lambda$
given by $\{f_i\}_{i=0}^2$ where
\bsenn
\left\{ \begin{array}{lll}
f_0(x) &=& \lambda x \\
f_1(x) &=& \lambda (x +1) \\
f_2(x) &=& \lambda (x+3)
\end{array} \right.
\esenn
Using techniques very different from ours, Pollicott and Simon \cite{PS}, in answer
to a question posed by Keane, have recently shown that for almost all $\lambda < 1/3$
the Hausdorff dimension $\Hdim({\cal G}_\lambda)$ is equal to $-\ln 3 /\ln \lambda$, while there is a dense subset of $[ 1/4, 1/3 ]$ such that if $\lambda$ belongs to this set then $\Hdim({\cal G}_\lambda)$ is strictly less than 
$-\ln 3 /\ln \lambda$.
(This is related to a problem posed by Erd\"os: for
which $\lambda \in [1/3,1]$ is the invariant density related to this system
singular with respect to Lebesgue measure (see \cite{So})).

For completeness, we list the definition of the Hausdorff dimension of a set here.
\begin{defn}
\label{def=Hdim}
Let ${\cal V}_\delta$ be the collection of covers of a set $C$ whose members
have diameter less than or equal to $\delta$. Let
\bsenn
{\cal H}^d_\delta(C) = \inf_{{\cal V}\in {\cal V}_\delta} \sum_{V_i \in {\cal V}} |V_i|^d \quad .
\esenn
The Hausdorff dimension $\Hdim (C)$ of $C$ is given by:
\bsenn
\Hdim(C) = \inf \{d \such \lim_{\delta \rightarrow 0} {\cal H}^d_\delta(C) =0\} \quad .
\esenn
\end{defn} \bigskip

\begin{defn}
\label{def=limcap} Let $C$ be a compact set in $\R^n$. The limit capacity (also known as the upper box counting dimension)
$d_c(C)$ of $C$ is given by:
\bsenn
d_c(C) = \limsup_{\delta \rightarrow 0} \dpy\frac{\ln \nu_C(\delta)}{-\ln \delta} \quad ,
\esenn
where $\nu_C(\delta)$ is the minimum number of balls of radius $\delta$ needed to cover 
C.
\end{defn}

\vskip.2in
\vskip .2in

\noindent
{\bf Acknowledgements:} We are grateful  to Marcelo Viana for useful
conversations on an early version of Theorem \ref{thm=cont-dim}.

\pagebreak

\section{Definitions}
\label{definitions}
\setcounter{figure}{0}
\setcounter{equation}{0}

\vskip.2in

Let $\bar{I}$ be the compact unit ball in $\R^n$. The space $H(\bar{I})$ is the 
collection of compact subsets of $\bar{I}$. For any compact $A$, denote 
its  $\epsilon$-neighborhood in $\bar{I}$ by $N_\epsilon (A)$. We topologize $H(\bar{I})$  by endowing it with  the so-called Hausdorff metric $\Hd$: 
\bsenn
\mbox{If}\, A_1, A_2 \in H(\bar{I}) \mbox{then} \Hd(A_1,A_2) = \max \{ \epsilon_1,\epsilon_2\} \where \epsilon_i = \inf \{ \epsilon | N_{\epsilon}(A_i) \supset A_j | j\neq i\} \quad .
\esenn
Observe that this metric is induced by 
the usual Euclidean distance $|\cdot |$ in $\R^n$. With this topology $H(\bar{I})$ 
becomes a complete, compact and metric space (see \cite{Hu}). 

The Hausdorff distance satisfies the following (strong) property:

\begin{lem}
\label{thm=hausdorff}
Let $A = \cup_{i} A_i$ and $B = \cup_{i} B_i$ be elements of $H(\bar{I})$. Then\\
i) $\Hd(A,B) \leq \sup_i \Hd (A_i,B_i)$ .\\
ii) $\Hd(A,B) \leq \sup_i \Hd (A,B_i)$ .
\end{lem}

\noindent
{\bf Proof:} Note that ii) follows immediately from i).
To prove i), note that for all $\epsilon$ greater than the right hand side, we have
\bsenn
A_i \subseteq N_\epsilon(B_i) \Rightarrow \cup_i A_i \;\subseteq \; \cup_i N_\epsilon (B_i) = 
N_\epsilon (\cup_i B_i),
\esenn
and vice versa.
\QED

\begin{defn}
\label{def=gifs}
A generalized iterated function system ${\cal F}$, is a map 
${\cal F}: H(\bar{I}) \rightarrow H(\bar{I})$ with the following properties: 
\begin{itemize}
\item  If $A \in H(\bar{I})$ then ${\cal F}(A) = \cup_{x \in A} {\cal F}(x)$ .
\item $\exists L \in [0,1)$ such that $\Hd({\cal F}(x), {\cal F}(y)) \leq L\cdot Hd(x,y)$ .
\end{itemize}
\end{defn}

\noindent
{\bf Remark:} One can convince oneself that not every generalized iterated 
function system is an iterated function system by considering the function 
 $z \rightarrow \lambda z^3$ on the unit disk, where $\lambda < 1/3$. It is still an open question whether 
every generalized iterated function system can be generated by continuous 
functions (as opposed to homeomorphisms) in the manner of Definition \ref{def=ifs}.

The last of the two properties in Definition \ref{def=gifs} will be referred to as contractiveness, 
since it implies that ${\cal F}$ is a contraction on $H(\bar{I})$. 

\begin{lem} 
\label{thm=contraction} A generalized iterated function system ${\cal F}$ 
is a contraction on $H({\bar I})$. In particular, ${\cal F}$ is continuous.
\end{lem}

\noindent
{\bf Proof:} Recall that $\epsilon$-neighborhoods in $\bar{I}$ are denoted 
by $N_\epsilon$. Suppose $B \subseteq 
N_\epsilon(A)$. By symmetry, it is sufficient to prove that
\bsenn
\Hd ({\cal F}(A),{\cal F}(\overline{N_\epsilon(A)}) \leq L\cdot \epsilon
\esenn
Now, using Definition \ref{def=gifs} and then applying Lemma \ref{thm=hausdorff} 
twice, one calculates:
\benn
\Hd ({\cal F}(A),{\cal F}(\overline{N_\epsilon(A)})) &\leq& \Hd (\cup_{a\in A}{\cal F}(a),\cup_{a\in A} {\cal F}(\overline{N_\epsilon(a)}) \leq \sup_{a \in A} \Hd ({\cal F}(a),{\cal F}(\overline{N_\epsilon(a)}) \\
&\leq& \sup_{a\in A} \sup_{b\in N_\epsilon (a)} \Hd ({\cal F}(a),{\cal F}(b)) \leq L\cdot \epsilon \quad .
\eenn
\QED

\begin{defn}
\label{def=G} The space $\G^0$ is the space of generalized iterated function systems
together with the following metric:
\bsenn
d_0({\cal F}, {\cal G}) = \max_{x \in {\bar I}} \Hd ({\cal F}(x),{\cal G}(x)) \quad .
\esenn
\end{defn}

$\G^0$ is complete and compact.  This is not difficult to show.  However, we will not use this fact, so we omit the proof. 

One usually defines an iterated function system (\cite{Fa}) 
as a system ${\cal F}$ consisting of a finite set of homeomorphisms. 
Even if one allows only  finite sets of affine transformations,  the theory of iterated function
systems is very rich and has many applications (see \cite{SB} for an overview).
Let $D_L^0(\bar{I})$ be space of contractions on ${\bar I}$ with Lipschitz constant 
$L$ and equipped with the usual sup-metric $\parallel \cdots \parallel _0$.

\begin{lem}
\label{thm=gifs}
If ${\cal F}$ is  a compact subset of $D_L^0(\bar{I})$ , then ${\cal F}$ satisfies
the hypotheses for a generalized iterated function system.
\end{lem}

\noindent
{\bf Proof:} We first prove that ${\cal F}: H(\bar{I}) \rightarrow H(\bar{I})$ 
is well-defined; that is, that the image of a compact set is compact.

To this end, define 
\bsenn
e : \bar{I} \times D^0_L(\bar{I}) \rightarrow H(\bar{I}) 
\esenn
by
\bsenn
e(x,f)=\{f(x)\} \quad .
\esenn
To prove that $e$ is continuous we let 
 $V$ be  open  and $(x, f) \in e^{-1}(V)$.  Then  $\{f(x)\} \in V$, and for some $\epsilon > 0$, $V$ contains an $(1+L)\epsilon$-neighbourhood of $\{f(x)\}$.  Now suppose $(z, g) \in N_{\epsilon}(x) \times N_{\epsilon}(f)$. Then 
\bsenn
|g(z) - f(x)| \leq |g(z) - f(z)| + | f(z) - f(x)| \\
\leq (1+L)\epsilon.
\esenn
That is, $(z, g) \in e^{-1}(V)$.
This shows that  $e$ is continuous.

Thus for  compact sets $A \subset \bar{I}$ and $\{f_j\}_{j \in J}\subset D^0_L(\bar{I})$, the set 
$\cup_{x\in A} \cup_{j\in J} f_j(x)$ is also compact and ${\cal F}$ is well-defined.

The remaining requirements are easy:
For $A \in H(\bar{I})$ it is true by definition that  ${\cal F}(A) = \cup_{x \in A} {\cal F}(x)$.  Furthermore, $\Hd({\cal F}(x), {\cal F}(y)) \leq \sup_{f \in {\cal F}}(f(x), f(y)) \leq L\cdot |x - y|$.  It now follows immediately that 
$\Hd({\cal F}(A), {\cal F}(B)) \leq L\cdot \Hd(A, B)$, and thus that ${\cal F}$ 
is continuous. 
\QED

It is not clear to us how to define, in any natural way, the notion of 
differentiability of a generalized iterated function system. 
So,  for our discussion of the Hausdorff dimension we consider only systems that consist of compact sets of diffeomorphisms (Definition \ref{def=G^1}).  Note, however, that the space we obtain is not complete.

\begin{defn}
\label{def=D-1-L} The metric space 
$D^{1+\alpha}_{L}({\bar I})$ is  the set of diffeomorphisms $f$ from 
${\bar I}$ to $f({\bar I}) \subseteq  I$, with Lipschitz constant $L < 1$, 
and with  $\alpha$-H\"older continuous derivatives. 
 We use the $C^1$-metric  on $D^{1+\alpha}_{L}({\bar I})$:  
\bsenn
\parallel f - g \parallel_1 = \sup_{x \in \overline{I}} \left( | f(x) - g(x) | + \parallel Df|_x - Dg|_x \parallel \right) \quad 
\esenn
\end{defn}

\begin{defn}
\label{def=difs}
A differentiable iterated function
system is a map ${\cal F} : H(\bar{I}) \rightarrow H(\bar{I})$, defined by 
\bsenn
{\cal F}(x) = \cup_{j \in J} f_j(x) \quad  ,
\esenn
where $\{f_j\}_{j\in J}$ a compact set in $D^{1+\alpha}_{L}(\bar{I})$. 
\end{defn}

  The $C^1$-metric on $D^{1+\alpha}_{L}({\bar I})$ induces a corresponding Hausdorff metric $\Hd^*$ on the space of compact subsets of $D^{1+\alpha}_{L}({\bar I})$.  In the following definition we use the same symbols ${\cal F}$ and ${\cal G}$ for  compact sets in $D^{1+\alpha}_{L}({\bar I})$ and for the differentiable iterated function systems they generate:

\begin{defn}
\label{def=G^1} The space $\G^1$ is the space of differentiable iterated function systems
together with the  metric  $\Hd^*({\cal F}, {\cal G})$ .
\end{defn}

\pagebreak

\section{The Measure Estimate}
\label{measure}
\setcounter{figure}{0}
\setcounter{equation}{0}
 
\vskip.2in

This section deals with the measure theoretic properties of the invariant 
sets. In particular, we prove that a generalized iterated function system 
${\cal F}$ is always a point of upper semi-continuity of the function that 
associates with ${\cal F}$ the Lebesgue measure of its invariant set. 
The modulus of semi-continuity is estimated in terms of the limit capacity 
of the boundary of the invariant set. (Counter-examples to continuity 
of the measure function are discussed in the introduction).

As is well-known, the fixed point is a continuous function of ${\cal F}$. More 
precisely:

\begin{prop}
\label{thm=cont} $\Hd (\Lambda({\cal F}), \Lambda({\cal G})) \leq \dpy \frac{d_0({\cal F}, {\cal G})}{1-L}$ .
\end{prop}

\noindent
{\bf Proof:} Observe that by the triangle inequality 
\benn
\Hd (\Lambda({\cal F}), \Lambda({\cal G})) &=& \Hd ({\cal F} \Lambda({\cal F}), {\cal G}\Lambda({\cal G})) \leq 
\Hd ({\cal F} \Lambda({\cal F}), {\cal F} \Lambda({\cal G})) +\Hd ({\cal F} \Lambda({\cal G}), {\cal G} \Lambda({\cal G}))   \\
&\leq& L\cdot \Hd (\Lambda({\cal F}), \Lambda({\cal G}))  + d_0({\cal F}, {\cal G}) \quad .
\eenn
\QED

\begin{theo}
\label{thm=cont_meas} The function $\mu$ is upper semi-continuous.
\end{theo}

\noindent
{\bf Proof:} The function $\mu$ is the composition of $\Lambda$ (which is continuous), 
and the Lebesgue measure function on $H({\bar I})$, which we will also denote by $\mu$ since no confusion is possible. It suffices to prove that the latter 
 is semi-continuous.

Suppose $\Lambda_0 \in H({\bar I})$ and suppose $\epsilon > 0$ is given.

Clearly, the neighborhoods $N_{1/n}(\Lambda_0)$ form a collection of monotone 
decreasing sets with 
\[
\lim_{n \rightarrow \infty} N_{1/n}(\Lambda_0) = \Lambda_0 \quad .
\]
 Since $\mu$ is a continuous measure this implies that 
\bsenn
\lim_{n\rightarrow\infty} \mu(N_{\frac{1}{n}})= \mu(\lim_{n\rightarrow\infty}(N_{\frac{1}{n}}))= \mu(\Lambda_0) \quad .
\esenn
Therefore, if n is large enough,
$\Hd (\Lambda, \Lambda_{0}) < \frac{1}{n}$, and so $\mu(\Lambda) < \mu(\Lambda_{0}) + \varepsilon$.
\QED

Notice, that the semi-continuity is not uniform (it can't be according to 
the examples in the introduction). But we can, in fact, 
estimate its modulus.

In the following,  $\partial \Lambda$ denotes the boundary of the set
$\Lambda$.

\begin{prop}
\label{thm=modulus} Suppose $\Lambda_{0} \in H(\bar{I})$, and $d > d_0 = d_c(\partial\Lambda_{0})$. 
Then for  any $\epsilon > 0$ and sufficiently small $\Delta > 0$ the following is true:
If $ \Lambda \in N_{\Delta}(\Lambda_{0})$, then
\bsenn
 \mu(\Lambda) \leq \mu(\Lambda_{0}) + \epsilon\Delta^{n-d} \quad .
\esenn 
\end{prop}

\noindent
{\bf Proof:}  
Observe that $\partial\Lambda_{0}$ can be covered by $p(\Delta)\Delta^{-d_0}$ 
balls of radius $\Delta$, where by definition 
\[ \limsup_{\Delta \rightarrow 0} \frac{\ln p(\Delta)}{-\ln \Delta} = 0 \quad .\] 
  If we increase the radius
of each of these balls to $3\Delta$, keeping their centers fixed, then the larger balls, together with $\Lambda_{0}$ will cover $N_{\Delta}(\Lambda_{0})$.  Thus, if $K_{n}$ is the volume of the unit ball in ${\R}^{n}$,
\benn
\mu(\Lambda)&\leq& \mu(\Lambda_{0}) + \mu(N_\Delta(\partial\Lambda_0))\\ 
&\leq& \mu(\Lambda_{0}) + (3\Delta)^{n}K_{n}\cdot p(\Delta)\Delta^{-d_0}.\\
&=& \mu(\Lambda_{0}) + (3\Delta)^{n}K_{n}\cdot (p(\Delta)\Delta^{\delta})
\cdot \Delta^{-d_0 -\delta} \quad .
\eenn
The term $(p(\Delta)\Delta^{\delta})$ tends to zero. Hence, for $\delta = d - d_0$ and $\Delta$ sufficiently small we have that the product $3^n \cdot  K_n \cdot (p(\Delta)\Delta^{\delta})$ is less than $\epsilon$.
\QED

\begin{cory}
Suppose ${\cal F}_0 \in \G^0$. Then for $d > d_c(\partial\Lambda({\cal F}_{0}))$, and for any $\epsilon > 0$ and  sufficiently small $\Delta > 0$ the following is true:
\[
d_0({\cal F}, {\cal F}_{0}) \leq \Delta \Rightarrow \mu(\Lambda({\cal F})) \leq \mu(\Lambda({\cal F}_{0})) + \epsilon\Delta^{n-d} \quad .
\]
\end{cory}

The corollary can be used to estimate the limit capacity in cases such as 
the example of Theorem \ref{thm=3branch}, 
where we know the measure of the invariant set but not its dimension.

\pagebreak

\section{The Distortion Estimate}
\label{distortion}
\setcounter{figure}{0}
\setcounter{equation}{0}

\vskip.2in

In one dimension, there is an elegant theory to obtain distortion estimates. 
This theory is described in various research papers and expository works 
(for example \cite{GJSTV}, \cite{MS} and \cite{Ve3}). The first step in this line of thought
is the following. Consider the forward orbit of an interval $I_0$ under 
a function $f$ and write $I_i = f(I_{i-1})$. 
Let $|\ln |Df| |_{\alpha}$ denote the $\alpha$-H\"older norm of the logarithm of the
derivative $Df$ of $f$ (restricted to the forward orbit of $I_0$). The distortion 
is the ratio of derivatives of high iterates of $f$ on a small interval $I_0$. 
It is given by the following expression (see the references listed above). 
\bse
\sup_{x_0,y_0 \in I_0} | \ln \left| \dpy\frac{Df^n(x_0)}{Df^n(y_0)} \right| | \leq 
|\ln |Df| |_{\alpha} \sum_i |I_i|^\alpha \quad .
\label{eq=distortion1}
\ese

In higher dimension, there is no convenient theory for the calculus 
of distortions. In this section, the beginnings of such a calculus 
are developed. Our theory mimics the derivation of the above estimate for the one-dimensional
calculus, but its elaboration is more awkward. Also we will iterate not necessarily 
by the same function, but each time by a function picked from an a priori
fixed (compact) set of functions. The latter generalization 
complicates the notation, but not the mathematics. 
We will need this theory in the next section where the dimension estimates 
are done.

Let us start by outlining the general idea of the estimates.
Suppose ${\cal F}$ is a compact set of diffeomorphisms in $D^1_{L, C}({\bar I})$.
Pick a ball $B_0$ of unit size and a sequence of contractions 
$f_i \in {\cal F}$. Define 
$f_i(B_{i-1}) = B_i$. We will express the distortion or nonlinearity as a sum of 
the diameters of the regions $B_i$, just as in the one-dimensional case,
except that now there will be a penalty for non-conformality. 
More precisely, choose two points $x_0$ and $y_0$ in $B_0$ and 
denote the images of $x_{i-1}$ and $y_{i-1}$ under $f_i$ by $x_i$ and $y_i$. The usual 
operator norm is written as $\parallel \cdot \parallel$. Define
\bse
C_n = C_n(f_n, \cdots f_1; x_0,y_0) = \left( D(f_n \cdots f_1)|_{x_0} \right)^{-1} \cdot \left( D(f_n \cdots f_1)|_{y_0} \right) \quad.
\label{eq=matrixC}
\ese
The idea is to obtain estimates for the logarithm of 
$\sup_{x_0,y_0 \in \bar{I},f_j \in {\cal F}} \parallel C_n \parallel$. 

Here is the main estimate of this section. To simplify the notation of its proof, 
we set for $n\geq 1$ 
\bse
\delta_n = (Df_{n}|_{x_{n-1}})^{-1} (Df_{n}|_{y_{n-1}}) - \Id  \quad .
\label{eq=delta}
\ese
We also agree that for a set $X$ in $\R^n$, we will denote its diameter by $|X|$. 
Finally, we simplify the notation for the conformality estimates.

\begin{defn}
\label{def=qn} Suppose ${\cal F} \in \G^1$ and define
\bsenn
Q(n) \equiv  \max_{i\leq n} \sup_{f_i \in {\cal F}, x_0 \in {\bar I}} \ln \parallel (Df_i \cdots f_1|_{x_0})\parallel \cdot \parallel (Df_i \cdots f_1|_{x_0})^{-1}\parallel \quad .
\esenn
\end{defn}

Note that by definition \ref{def=semi-conformal}, ${\cal F}$ is called 
semi-conformal if $Q(n)/n$ tends to zero.

\begin{theo}
\label{thm=distortion}
Let ${\cal F}\in \G^1$.  
Then  there is a constant $C$ such that for all vectors $v \in S^{n-1}$
\bsenn
|\ln |C_n v|| \leq C \sum_{j=0}^{n-1} e^{Q(j)} \cdot |B_{j-1}|^\alpha \quad .
\esenn
\end{theo}

\noindent
{\bf Proof:} Observe that 
\benn
C_{n+1} &=& ( D(f_n \cdots f_1)|_{x_0} )^{-1} \cdot (\Id + \delta_{n+1}) \cdot ( D(f_n \cdots f_1)|_{y_0} ) \\
&=& ( D(f_n \cdots f_1)|_{x_0} )^{-1} \cdot (\Id + \delta_{n+1}) \cdot ( D(f_n \cdots f_1)|_{x_0} ) C_n \\
&=& ( \Id + (D(f_n \cdots f_1)|_{x_0} )^{-1} \cdot \delta_{n+1}  \cdot D(f_n \cdots f_1)|_{x_0})  C_n \\
&=& ( \Id + (D(f_n \cdots f_1)|_{x_0} )^{-1} \cdot \delta_{n+1}  \cdot D(f_n \cdots f_1)|_{x_0}) \\
&& \quad \cdot ( \Id + (D(f_{n-1} \cdots f_1)|_{x_0} )^{-1} \cdot \delta_{n} \cdot  D(f_{n-1} \cdots f_1)|_{x_0})\cdots (\Id + \delta_1) \quad .
\eenn
Further, because ${\cal F}$ and $\bar{I}$ are compact, there is a
uniform upper estimate for $\|(Df|_x)^{-1}\|$ on ${\cal F} \times \bar{I}$.
By the assumption  of uniform H\"{o}lder continuity it then follows that 
\bsenn
\parallel \delta_i\parallel \leq C |B_{i-1}|^\alpha \quad .
\esenn
Now, using Schwarz' inequality (twice) and the triangle inequality, we obtain that
\bse
\parallel C_{n} \parallel \leq \dpy\prod_{j=1}^n \left( 1 + \parallel (Df_j \cdots f_1|_{x_0})^{-1}\parallel \cdot \parallel \delta_{j+1} \parallel \cdot \parallel Df_j \cdots f_1|_{x_0}\parallel \right) \quad . 
\label{eq=nonlinearity} 
\ese

Note that the matrices $C_n$ are invertible and that
\bsenn
(C_n(f_n, \cdots f_1; x_0,y_0))^{-1} = C_n(f_n, \cdots f_1; y_0,x_0) \quad .
\esenn
Thus $\parallel C_n \parallel^{-1}$ also satisfies equation (\ref{eq=nonlinearity}). 
Since 
\bsenn
\parallel C_n^{-1} \parallel^{-1} \leq |C_nv| \leq \parallel C_n \parallel \quad ,
\esenn
we now obtain the estimate for $|\ln |C_nv||$ upon taking logarithms. 
\QED 

\noindent
{\bf Remark:} One will notice that 
the estimate also holds for expanding diffeomorphisms. The 
problem is that in this case the estimate is bad, since the factor multiplying 
the largest of the diameters is uncontrollable. It is, however, possible to
derive a theorem along the same lines for expanding diffeomorphisms. This 
is done by redefining $C_n$ and $\delta_n$ as follows:
\benn
C_n &=& \left( D(f_n \cdots f_1)|_{x_0} \right) \cdot \left( D(f_n \cdots f_1)|_{y_0} \right)^{-1} \quad.\\
\logand \quad \delta_n &=& (Df_n|_{x_{n-1}})\cdot (Df_n|_{y_{n-1}})^{-1} - \Id \quad .
\eenn
(Notice that in both cases we write the expanding term first.) Now one derives 
from the recursion:
\bsenn
C_{n+1} = (\Id + \delta_{n+1}) \cdot (Df_{n+1}|_{y_n})^{-1}\cdot C_n \cdot (Df_{n+1}|_{y_n}) \quad ,
\esenn
and continues the reasoning as in the above proof.

\vskip.2in
\noindent
{\bf Remark:} Notice that we assume that derivatives and their inverses exist,
precluding singularities. The one-dimensional theory also includes a calculus of distortions
in the presence of singularities in the derivative by using for example
Koebe estimates or cross-ratios (reviewed in \cite{MS}). 
At present there are no tools to study
the analogous problem in higher dimension.

\vskip.2in
\noindent
{\bf Remark:} The term $e^{Q(j)}$ in the theorem 
has value one if the $f$'s are conformal. In reality this factor is a penalty 
for deviation from conformality. It is easy to see that if the moduli of the 
eigenvalues of $D(f_n\cdots f_1)$ bunch together sufficiently, then the lack 
of conformality will be compensated for  by the exponential decrease of the $|B_i|$. 
In particular, it is clear that a semi-conformal system has bounded distortion. 

\vskip .2in

\begin{defn}
\label{def=bdd-dist} For ${\cal F} \in \G^1$, define the distortion as
\bsenn
D(n) = \max_{i\leq n} \sup_{x,y \in {\bar I}, f_j \in {\cal F}} |\ln \| C_i(f_i, \cdots f_1; x,y)\| | \quad .
\esenn
The system is said to have bounded distortion if $D(n)/n$ tends to zero.
\end{defn}

Note that this is slightly more general than the usual requirement which is that
$D(n)$ is uniformly bounded. 

\begin{cory}
\label{thm=bound-distortion}
Let ${\cal F} \in \G^1$.  Then for all $x_0$ and $y_0$ in ${\bar I}$ 
and each unit-vector $v$ 
there is a unit-vector $w(v)$ such that  
\bsenn
\left| \ln \dpy\frac{|D(f_n \cdots f_1)|_{x_0} (w(v))|}{|D(f_n \cdots f_1)|_{y_0} (v)|} \right| \leq D(n) \quad .
\esenn
\end{cory}

\noindent
{\bf Proof:} Note that 
\bsenn
D(f_n\cdots f_1)|_{x_0} (C_nv) = D(f_n\cdots f_1)|_{y_0} (v)
\esenn
We choose $w(v) = \frac{C_nv}{|C_nv|}$. Then 
\bsenn
\dpy\frac{|D(f_n \cdots f_1)|_{x_0} (w(v))|}{|D(f_n \cdots f_1)|_{y_0} (v)|} = 
\dpy\frac{1}{|C_nv|} \quad .
\esenn
Now take logarithms and apply the theorem.
\QED

To obtain good bounds in the dimension calculations, we need to know 
that the asymptotic rate of contraction is independent 
of the direction in the tangent space. This is implied by semi-conformality 
(see definition \ref{def=semi-conformal}) and the following result:

\begin{cory}
\label{cor=semi-conformal} 
Let ${\cal F} \in \G^1$. Then for all unit-vectors $v$ and $u$:
\bsenn
\left| \ln \dpy\frac{|D(f_n \cdots f_1)|_{x_0}(u)|}{|D(f_n \cdots f_1)|_{y_0}(v)|}\right|
\leq D(n) + Q(n) \quad .
\esenn
\end{cory}

\noindent
{\bf Proof:} We write this as 
\bsenn
\left| \ln \dpy\frac{|D(f_n \cdots f_1)|_{x_0} (w(v))|}{|D(f_n \cdots f_1)|_{y_0} (v)|}\right| +
\ln \dpy\left|\frac{|D(f_n \cdots f_1)|_{y_0} (u)|}{|D(f_n \cdots f_1)|_{y_0} (w(v))|} \right| \quad .
\esenn
To the first term we apply the previous result. The second is calculated
with the help of
\bsenn
\parallel (Df_n \cdots f_1|_{x_0})^{-1} \parallel^{-1} |v| \leq |(Df_n \cdots f_1|_{x_0})v| \leq \parallel Df_n \cdots f_1|_{x_0}\parallel |v| \quad ,
\esenn
and the definition of $Q(n)$. 
\QED

In the remainder of this section, we discuss the relation between derivatives
and sizes of domains. In particular, we derive a version of the mean value theorem to 
to obtain better estimates of the diameter of the iterate of a region.

\begin{lem}
Let $A$ and $B$ be connected compact sets in $\R^n$ and suppose in addition that 
$A$ is convex.
Suppose that $g: A \rightarrow B$ is a diffeomorphism. Then there
is a point $a_+ \in A$ and a $v_{a_+} \in T_{a_+}A$ (the tangent space of $A$ at $a_+$) 
such that
\bsenn
\dpy\frac{|Dg|_{a_+} v_{a_+}|}{|v_{a_+}|} \geq \dpy\frac{|B|}{|A|} \quad .
\esenn
\end{lem}

\noindent
{\bf Proof:} Let $w$ and $z$ in $B$ such that $|w - z| = |B|$ and
let $x = g^{-1}(w)$ and $y = g^{-1}(z)$. Connect $x$ and $y$ 
by a straight segment $\gamma \in A$ (by the convexity of $A$) and  
parametrize this curve by arclength ($|D\gamma| = 1$). 
Then 
\bsenn
|B| = \left| \dpy\int_0^{|x - y|} Dg(\gamma(t))\cdot D\gamma(t) \, dt \right| \leq |A| \cdot \max_{x\in A} \parallel Dg_x\parallel \quad .
\esenn
Now choose $a_+$ to be the point where the maximum is assumed. 
\QED

\begin{lem}
Let $A$ be a closed ball and $B$ be a set in $\R^n$.
Suppose that $g: A \rightarrow B$ is a diffeomorphism. Then there 
is a point $a_- \in A$ and a $v_{a_-} \in T_{a_-}A$ such that 
\bsenn
\dpy\frac{|Dg_{a_-} v_{a_-}|}{|v_{a_-}|} \leq \dpy\frac{|B|}{|A|} \quad . 
\esenn    
\end{lem} 
 
\noindent 
{\bf Proof:} From elementary calculus, we know that
\bsenn
\dpy\int_A \dpy\frac{|\det Dg_x|}{{\rm vol}(A)} \, d^n x = \dpy\frac{{\rm vol}(B)}{{\rm vol}(A)} \quad .
\esenn
The right hand side of this equation is the average of the positive function 
$|\det Dg_x|$. Thus there is a $a_- \in A$ such that
\bsenn
|\det Dg|_{a_-}| \leq \dpy\frac{{\rm vol}(B)}{{\rm vol}(A)} \quad .
\esenn
Denote the eigenvalues of $Dg|_{a_-}$ by $\{\lambda_i\}_{i=1}^{n}$ (counting
multiplicity). Observe that ${\rm vol}(B)$ is no greater than the volume of a 
ball with diameter $|B|$. Thus
\bsenn
\prod_{i=1}^{n}|\lambda_i| \leq \dpy\frac{|B|^n}{|A|^n} \quad .
\esenn
By taking logarithms and dividing by $n$, it becomes obvious that the average 
of $\{\ln|\lambda_i|\}_{i=1}^{n}$ is no greater than $\frac{|B|}{|A|}$. Thus 
there must be an eigenspace $V_{a_-}$ of $Dg_{a_-}$ satisfying the lemma.
\QED

These two lemmas imply the higher dimensional version of the mean value theorem that 
we will use in the next section.

\begin{cory}
\label{thm=rolle}
Let $A$ be a closed ball and $B$ a set in $\R^n$. 
Suppose that $g: A \rightarrow B$ is a diffeomorphism. Then there  
is a point $a \in A$ and a $v \in T_aA$ such that   
\bsenn 
\dpy\frac{|Dg|_a v|}{|v|} = \dpy\frac{|B|}{|A|} \quad .  
\esenn     
\end{cory}

\noindent
{\bf Proof:} Note that the transformation
\bsenn
T:  A \times TA \rightarrow \R^+
\esenn
defined by
\bsenn
T(x,v) = \dpy\frac{|Dg|_x v|}{|v|}
\esenn
is continuous and $A \times TA$ is path-connected. The result is thus a consequence 
of the previous lemmas.
\QED

We now use these results to derive a general statement about scalings 
in contracting maps. 

\begin{theo}
\label{thm=scaling}
 Let ${\cal F} \in \G^1$ and suppose that $f_a$ is a composition of at most $n$ functions of ${\cal F}$  and $f_b$ is a composition of arbitrary length.  Then for any ball $B$ we have 
\bsenn
\left| \ln \left( \dpy\frac{|f_a f_b (B)|}{|f_a(B)|} \cdot \frac{|B|}{|f_b(B)|} \right) \right|
 \leq 2Q(n) + 2D(n) \quad.
\esenn
\end{theo}

\noindent
{\bf Proof:}  The expression in the theorem can be written as:
\bsenn
\left| \ln \left( \dpy\frac{|f_a f_b (B)|}{|B|} \cdot \frac{|B|}{|f_a(B)|} \cdot \frac{|B|}{|f_b(B)|} \right) \right|
\esenn
With the help of Corollary \ref{thm=rolle}, we get
\bsenn
\frac{|f_a f_b (B)|}{|B|} = \frac{|(Df_a \cdot Df_b)|_x\ v_x|}{|v_x|} =
\frac{|Df_a|_yv_y|}{|v_y|}\ \frac{|Df_b|_xv_x|}{|v_x|} \quad ,
\esenn
where $v_y$ is a unit vector in the direction of $Df_b|_xv_x$.  The other derivatives can also be calculated with the help of the same corollary, to give 
\bsenn
\left( \dpy\frac{|f_a f_b (B)|}{|f_a(B)|} \cdot \frac{|B|}{|f_b(B)|} \right)
= \frac {| Df_a|_{x_1}(v_1)|}{|Df_a|_{x_2}(v_2)|} \cdot  \frac                 {|Df_b|_{y_1}(w_1)|}{|Df_b|_{y_2}(w_2)|} \quad .
\esenn
The result now follows from Corollary \ref{cor=semi-conformal}.
\QED

\vskip.2in

\pagebreak

\section{The Dimension Estimate}
\label{dimension}
\setcounter{figure}{0}
\setcounter{equation}{0}

\vskip.2in

We prove that if ${\cal F} \in \G^1$ is a semi-conformal differentiable iterated function system, then it is a point of lower semi-continuity of the function 
that evaluates the Hausdorff dimension. 
We note here that if ${\cal F}$ is a one-dimensional system with finitely 
many branches and satisfying a strong condition on the distance of the individual 
branches, then the Hausdorff dimension 
varies continuously (see \cite{Ta} or \cite{Ru} for more information). Without that condition,
it is clear that the dimension is not continuous as observed in the introduction. 
Nonetheless, the proof of the semi-continuity given here bears resemblance to Takens' proof. 

For a given system ${\cal F} \in \G^1$, we choose positive constants $K = K({\cal F})$ and $k = k({\cal F})$ such that for all $f \in {\cal F}$ and $x \in \overline{I}$\ ,
\bse
\begin{array}{ll}
 \left(  \parallel (Df|_x)^{-1} \parallel\right)^{-1} & > e^{-K}\quad ; \\
 \parallel Df|_x \parallel & < e^{-k}\quad . 
\end{array}
\label{eq=constants}
\ese
Note that by continuity (\ref{eq=constants}) automatically holds for all  ${\cal F^{\prime}}$ in a $\G^1$-neighbourhood of ${\cal F}$ (see Definition \ref{def=G^1}). 

A dynamic cover ${\cal U}$ of $\Lambda({\cal F})$ is a finite cover by open sets each of 
which can be written as $f_n\cdots f_1(I)$. 

\begin{lem}
\label{thm=cover} Suppose ${\cal F} \in \G^1$.  Then for each $n > 0$, 
there is a dynamic cover ${\cal U}_n$ of $\Lambda({\cal F})$ such that for all
$U \in {\cal U}_n$: 
\bsenn
2 e^{-2Q(n)-2D(n)-K-nk} \leq |U| < 2 e^{-nk} \quad .
\esenn
Furthermore, all elements of ${\cal U}_n$ are of the form $f_m\cdots f_1(I)$ 
with $m \leq n$.
\end{lem}

\noindent
{\bf Proof:}  Since $\Lambda({\cal F}) \subset I$, \, $\{f(I)\}_{f \in {\cal F}} $ covers $\Lambda({\cal F})$.  Let ${\cal U}_1$ be a finite subcover.  It clearly satisfies the lemma. 
Denote the finite set of contractions selected in this process by ${\cal F}_1$. Thus
\bsenn
\Lambda({\cal F}) \subset \{f(I)\}_{f\in {\cal F}_1} = {\cal F}_1(I) \quad .
\esenn
We now continue by induction. Suppose  that for $i \in \{1,\cdots,n\}$, we have 
finite covers ${\cal U}_i = {\cal F}_i(I)$ of $\Lambda({\cal F})$ satisfying the lemma.
To construct ${\cal U}_{n+1}$, we will for each $W = f_w(I) = f_{\ell} \cdots f_{1}(I) \in {\cal U}_n$ find a finite cover of $W \cap \Lambda({\cal F})$ by sets of the form $f_wf_\alpha(I)$ where each $f_\alpha$ is also a composition of elements of ${\cal F}$.
 
If $W = \Lambda({\cal F})$ already satisfies the inequalities of the lemma for level $n+1$, we accept $W$ itself as a member of ${\cal U}_{n+1}$.  If this is not the case, we have
\bse
2e^{-(n+1)k} \leq |W| < 2e^{-nk}.
\label{eq=pinch}
\ese
Note that for any set $W = f_w(I)$ satisfying (\ref{eq=pinch}) we necessarily have $|w| \leq n$, where $|w|$ denotes the length of the composition $f_w$.  

We now replace $f_w(I)$ by $\{f_wf(I)\}_{f \in {\cal F}_{1}}$.  The latter is clearly a covering of $W \cap \Lambda({\cal F})$.  For some $ m \leq n$ we have, by  Theorem \ref{thm=scaling}, the following situation:
\bse
\begin{array}{lll}
|f_wf(I)| &=& C_w\dpy\frac{|f(I)|}{|I|} |f_w(I)| \quad,\\[.3cm]
|\ln C_w| &\leq & 2Q(m) + 2D(m) \quad ,\\[.3cm]
e^{-K}\leq &\dpy\frac{|f(I)|}{|I|} & < e^{-k} \quad . 
\end{array}
\label{eq=diameterrecursion}
\ese
Combining this with (\ref{eq=pinch}), we obtain
\bse
2e^{-(n+1)k - K -2Q(m)-2D(m)} \leq |f_wf(I)| < 2 e^{-(n+1)k +2Q(m)+2D(m)} \quad .
\label{eq=diameterrecursion2}
\ese
In particular, we have the necessary lower bound for  $ |f_wf(I)|$.
If the upper bound is  also satisfied we accept $f_wf(I)$ as a member of ${\cal U}_{n+1}$.   If not, then $f_wf(I)$ also satisfies (\ref{eq=pinch}).  This means that we can proceed, as for $f_w(I)$, by replacing $f_wf(I)$ by the covering
$\{f_wfg(I)\}_{g \in {\cal F}_1}$ of $f_wf(I) \cap \Lambda({\cal F})$, and treating its members as before.

Since
\bsenn
|f_wf_\alpha(I)| < 2 e^{-(n+1)k} \ \mbox{if}\ |w| + |\alpha| = n+1 \quad , 
\esenn
 this iterative procedure must terminate in success after no more than $n+1-|w|$ steps, resulting in a dynamic cover ${\cal U}_{n+1}$ satisfying the lemma.
\QED

With each dynamic cover ${\cal U}_n = {\cal F}_n(I)$ of $\Lambda({\cal F})$, associate a compact 
subset $\Lambda_n$ of $\Lambda$ in the following way: Let ${\cal V}_n$ be a 
maximal collection of disjoint members of ${\cal U}_n$. Thus by construction 
each $V\in {\cal V}_n$ is of the form $f(I)$ where $f$ belongs to a subset 
${\cal G}_n$ of ${\cal F}_n$. Clearly, ${\cal G}_n$ is a differentiable 
iterated function system consisting of a finite number, say $N_n$, of 
diffeomorphisms, and its invariant set $\Lambda_n=\Lambda({\cal G}_n)$ 
is a subset of $\Lambda$. 
The main reason for introducing the systems ${\cal G}_n$, is, of course, 
that they satisfy the open set condition (definition \ref{def=open-set}).
Therefore, the Hausdorff dimensions of their invariant sets are easy to calculate.
These dimensions serve as approximations to the dimension of $\Lambda$ (see
Proposition \ref{thm=dimension-estimate}). 

\begin{lem}
For each $V$ in ${\cal V}_n$, choose $x_V \in V$ and let ${\tilde V}$ be the ball with center $x_V$ but with radius $4e^{-nk}$. Then the collection 
$\tilde{{\cal V}}_n$ of these sets covers $\Lambda$, and for any $V \in {\cal V}_n$ we have 
\[ | {\tilde V} | \leq 4 e^{2Q(n) + 2D(n) + K} | V | \quad . \]
\end{lem}

The proof is easy and is therefore left to the reader. 

\begin{lem}
\label{thm=dim-est1}
Suppose ${\cal F} \in \S \cap \G^1$ and let ${\cal G}_n$ be the differentiable iterated function system constructed as above,
and suppose it consists of $N_n$ diffeomorphisms. Then
\bsenn
\Hdim (\Lambda_n) \leq \Hdim(\Lambda) \leq \limsup_{n\rightarrow \infty} \frac{\ln N_n}{nk} \quad .
\esenn
\end{lem}

\noindent
{\bf Proof:} The first of the inequalities is obvious because of the inclusion 
$\Lambda_n \subseteq \Lambda$.  For the second inequality let $d_n = \frac{\ln N_n}{nk}$.  Then for any $\epsilon > 0$,

\bse
\begin{array}{lllll}
{\cal H}^{\epsilon + d_n}_{8e^{-nk}}(\Lambda) &\leq& \sum_{\tilde{V} \in \tilde{{\cal V}_n}} | \tilde{V}|^{\epsilon + d_n} &\leq& \sum_{V \in {\cal V}_n} 4^{\epsilon + d_n} e^{(K + 2Q(n) + 2D(n))(\epsilon + d_n)} |V|^{\epsilon + d_n} \\[0.4cm]
& & &\leq& 8^{\epsilon + d_n} N_n e^{(K + 2Q(n) + 2D(n) - nk)(\epsilon + d_n)} \\[0.4cm]
& & &\leq& 8^{\epsilon + d_n} e^{K + 2Q(n) + 2D(n))(\epsilon + d_n) - nk\epsilon} \quad . 
\end{array}
\ese
This tends to zero as $n$ goes to infinity, and so establishes the upper bound for the Hausdorff dimension.
\QED

The actual calculation of the Hausdorff dimension uses the following 
result (see \cite{Fa}). 

\begin{prop}
\label{thm=fa}
Let the system ${\cal H} \in \G^1$ be a set of $N$ contractions 
satisfying the open set condition. Suppose further that
\bsenn
0 < e^{\lambda_-} \leq \dpy\frac{|Dh\cdot v|}{|v|} \leq e^{\lambda_+} < 1 \quad ,
\esenn
for all $h \in {\cal H}$. Then we have 
\bsenn
\dpy\frac{-\ln N}{\lambda_-} \leq \Hdim(\Lambda({\cal H})) \leq \dpy\frac{-\ln N}{\lambda_+} \quad .
\esenn
\end{prop}

\begin{prop}
\label{thm=dimension-estimate}
Consider the invariant sets $\Lambda_n$ of the systems ${\cal G}_n$ derived 
from the system ${\cal F} \in \S \cap \G^1$ and consisting of
$N_n$ contractions satisfying the open set condition. Let $k$ be the constant 
defined in equation (\ref{eq=constants}). Then we have: 
\bsenn
\lim_{n \rightarrow \infty} |\Hdim(\Lambda_n) - \dpy\frac{\ln N_n}{nk}| = 0 \quad .
\esenn 
\end{prop}

\noindent
{\bf Proof:} The proof consists of estimating the numbers $\lambda_-$ and 
$\lambda_+$ of the previous proposition for the sets $\Lambda({\cal G}_n)$. 
Let $g: I \rightarrow V \in {\cal V}_n$ be a member of the finite family 
${\cal G}_n$. By construction, the map $g$ is a composition of $m \leq n$ diffeomorphisms
$f \in {\cal F}$. For any $x \in {\bar I}$ and $v$ in the tangent space $T_x{\bar I}$, 
we split the basic estimate of Proposition \ref{thm=fa} into two parts:
\bsenn
\dpy\frac{|Dg_x\cdot v|}{|v|} = \dpy\frac{|Dg_a\cdot v_a|}{|v_a|}\cdot 
\dpy\frac{|Dg_x\cdot v|}{|Dg_a\cdot v_a|} \quad ,
\esenn
Where  $a$ and the unit tangent vector $v_a$ are selected so that 
the first ratio on the right hand side is estimated by applying the mean value theorem
(Corollary \ref{thm=rolle}). The second ratio is estimated by using Corollary \ref{cor=semi-conformal}.   Using Lemma \ref{thm=cover} these estimates give us
\bsenn
- nk - K - 3D(n) - 3Q(n) \leq \ln \dpy\frac{|Dg_x\cdot v|}{|v|} \leq  - nk + D(n) +Q(n) \quad .
\label{eq=dim-est}
\esenn

thus, by Proposition \ref{thm=fa}, 
\bse
\frac{\ln N_n}{nk + K + 3D(n) + 3Q(n)} \leq \Hdim(\Lambda({\cal G}_n)) \leq \frac{\ln N_n}{nk - D(n) -Q(n)} \quad .
\label{eq=hdimn}
\ese
Using semi-conformality, as $n \rightarrow \infty$,
this establishes the result. 
\QED

Recall the definition of the limit capacity, given in Section 1.  The limit capacity is always at least as big as the Hausdorff
dimension, because for the former we insist that the covering sets all have the
same diameter.  The following is an immediate consequence of Proposition \ref{thm=dimension-estimate} and Lemma \ref{thm=dim-est1}:

\begin{cory}
\label{thm=lim_cap}
Suppose that ${\cal F} \in \S\, \cap\G^1$.  Then the limit capacity of $\Lambda({\cal F})$ is equal to its Hausdorff dimension.
\end{cory}

In the mean time we have everything in place to prove the extension of the continuity 
result (one part of Theorem \ref{thm=main2}). The methods are exactly the same 
as those used in the previous proposition. 

\begin{theo}
\label{thm=main2-1}
Every point of  $\G^1 \cap \O \cap \S$ is a point of continuity of the function $\Hdim$ on $\G^1 \cap \O$.
\end{theo}

\noindent
{\bf Proof:} For ${\cal F} \in \G^1 \cap \O \cap \S$ and  ${\cal F}^{\prime} \in \G^1 \cap \O$ we let 
${\cal F}_n$ and ${\cal F}^{\prime}_n$ be the collection of iterates associated with the 
dynamic covers of $\Lambda({\cal F})$ and $\Lambda({\cal F}^{\prime})$ respectively.  If ${\cal F}$ and ${\cal F}^{\prime}$ are sufficiently close in the $\G^1$-topology, then they will satisfy (\ref{eq=constants}) for the same constants $k$ and $K$. 
  Because ${\cal F}, {\cal F}^{\prime} \in  \O$ we have $\Lambda({\cal F}_n) = 
\Lambda({\cal F})$   and $\Lambda({\cal F}^{\prime}_n) = \Lambda({\cal F}^{\prime})$.  We let 
$N_n = \#({\cal F}_n)$, and we define $D(n)$ and $Q(n)$ for ${\cal F}$ as in  Definitions \ref{def=bdd-dist} and \ref{def=qn}.  For
 ${\cal F}^{\prime}$ we denote the corresponding quantities by  
$D^{\prime}(n)$ and $Q^{\prime}(n)$. Then by (\ref{eq=hdimn}) in the proof of 
Proposition \ref{thm=dimension-estimate},
\bsenn
\frac{\ln N_n}{nk + K +3D(n)+3Q(n)} \leq \Hdim(\Lambda({\cal F})) \leq 
\frac{\ln N_n}{nk - D(n) - Q(n)}
\esenn
for each $n$.  A similar inequality holds for $\Hdim(\Lambda({\cal F}^{\prime}))$.

Fix $n$.  Then for ${\cal F}^{\prime}$ in a sufficiently small neighborhood ${\cal N}$ of ${\cal F}$ in \ $\G^0 \cap \O \cap \S$ \ it follows that $\#({\cal F}^{\prime}_n) = N_n$, and that  $D^{\prime}(n)$ and $Q^{\prime}(n)$ are arbitrarily close to $D(n)$ and $Q(n)$ respectively.  Thus, by choosing $n$ large initially, and ${\cal N}$ small, we can make 
$| \Hdim(\Lambda({\cal F}) - \Hdim (\Lambda({\cal F}^{\prime}))|$ as 
small as desired.
\QED

Here is the remaining half of the main result.

\begin{theo}
\label{thm=cont-dim}
Every point of  $\G^1 \cap \S$ is a point of lower semi-continuity of the function $\Hdim$ on $\G^1$.
\end{theo}

\noindent
{\bf Proof:} Recall that $\Lambda_n = \Lambda_n({\cal F})$ is the invariant set of the system ${\cal G}_n = {\cal G}_n({\cal F})$ 
satisfying the open set condition.  For a given ${\cal F} \in \G^1 \cap \S$  choose $n$ so that 
\bsenn
| \Hdim(\Lambda({\cal F})) - \Hdim(\Lambda_n({\cal F})) | 
< \epsilon/2 \quad .
\esenn
This is possible by Proposition \ref{thm=dimension-estimate}.
If we choose ${\cal F}^{\prime} \in \G^1$ sufficiently close to ${\cal F}$, then ${\cal G}_n({\cal F}^{\prime})$ is close to  ${\cal G}_n({\cal F})$ and  ${\cal G}_n({\cal F}^{\prime})$ also satisfies the open set condition.  We assume, using Theorem \ref{thm=main2-1}, that $|\Hdim(\Lambda_n({\cal F}^{\prime})) - \Hdim(\Lambda_n({\cal F}))| < \epsilon/2$.  Then, using Lemma \ref{thm=dim-est1},
\benn
\Hdim(\Lambda({\cal F}^{\prime})) &\geq& \Hdim(\Lambda_n({\cal F}^{\prime}))  \\ &\geq& \Hdim(\Lambda_n({\cal F})) - \epsilon/2 \\ &\geq& \Hdim(\Lambda({\cal F})) - \epsilon \quad . \\
\eenn
\QED

\noindent
{\bf Remark:} We note that in the proof of this theorem the imposition of the most 
restrictive condition, namely semi-conformality, arises  not 
from our limited knowledge of the calculus of distortions in higher 
dimension, but from the absence of methods for calculating the
Hausdorff dimension for sets in $\R^n$ with $n>1$. Admittedly there are
some methods that apply to certain affinely generated sets (see \cite{McM}) but these cannot be used here. 
It appears likely
that the semi-continuity of the theorem holds in a wider context; that is,
that the requirement of semi-conformality can be relaxed or dropped 
altogether. To address this problem, a deeper knowledge of the function
$\Hdim$ is required. Conceivably, this type of result is more readily 
proved using a different definition of the dimension.

\vfil \eject

\end{document}